\documentclass[notitlepage,11pt]{article}
\usepackage{amssymb,amsmath,comment}
\catcode`\@=11
\@addtoreset{equation}{section}

\catcode`\@=12

\usepackage{graphicx}

 \newtheorem{thm}{Theorem}[section]
 \newtheorem{cor}[thm]{Corollary}
 \newtheorem{lem}[thm]{Lemma}
 
\usepackage{comment}

\def\R{\mathbb{R}}

\def\S{\mathbb{S}}

\def\irn{\int_{\R^n}}

\def\da{\delta_{n,\alpha}}
\def\darad{\delta_{n,\alpha}^{\rm rad}}
\def\ha{h_{n,\alpha}}

\def\ga{\gamma_{n,\alpha}}

\def\Ma{M_{n,\alpha}}

\def\proof{\noindent{\textbf{Proof. }}}
\def\qed{\hfill {$\square$}\goodbreak \medskip}

\begin{document}

\title
{A class of second order\\
dilation invariant inequalities}

\author{Paolo Caldiroli\footnote{Dipartimento di Matematica, Universit\`a di Torino, via Carlo Alberto, 10 -- 10123 Torino, Italy, Email: {paolo.caldiroli@unito.it}.
Partially supported by the PRIN2009 grant "Critical Point Theory and Perturbative Methods for Nonlinear Differential Equations".} and
{Roberta Musina\footnote{Dipartimento di Matematica ed Informatica, Universit\`a di Udine, via delle Scienze, 206 -- 33100 Udine, Italy. Email: {roberta.musina@uniud.it}.
Partially supported by Miur-PRIN 2009WRJ3W7-001 ``Fenomeni di concentrazione e {pro\-ble\-mi} di analisi geometrica''}}}

\date{}

\maketitle

\begin{abstract}
\footnotesize We compute the best constants in some dilation invariant inequalities
for the weighted
$L^2$-norms of $-\Delta u$ and $\nabla u$, with weights being
powers of the distance from the origin.

\medskip

\medskip

\end{abstract}

\bigskip\bigskip\bigskip


\section{Introduction}
In recent years, there has been a growing interest in 
dilation invariant inequalities that are somehow related with the famous Rellich inequality
 \cite{Rel54}, \cite{Rel69}.
We shall not attempt to provide a complete list
of references on this subject. However, among the more recent
contributions we cite \cite{AGS}--\cite{Mor}, \cite{TerZog}
and references therein.

In the present paper we study a class of inequalities for the weighted $L^2$-norms of $-\Delta u$ and $\nabla u$.
More precisely, let $n\ge 2$ be a given integer, let $\alpha\in\R$ be a varying parameter, and let $\Sigma$ be a regular domain in $\S^{n-1}$. We are interested in inequalities of the form
\begin{equation}
\label{eq:DN}
\int_{\mathcal C_{\Sigma}}|x|^{\alpha}|\Delta u|^{2}dx\ge c
\int_{\mathcal C_{\Sigma}}|x|^{\alpha-2}|\nabla u|^{2}dx\quad\textrm{for any
$u\in C^2_c(\overline{\mathcal C_{\Sigma}}\setminus\{0\})$}
\end{equation}
where $\mathcal C_{\Sigma}$ denotes the cone in $\R^{n}$ spanned by $\Sigma$, namely
$$
\mathcal C_{\Sigma}=\left\{x\in\R^{n}\setminus\{0\}~\left|~\frac{x}{|x|}\in\Sigma\right.\right\}.
$$
Notice that $\mathcal C_{\Sigma}=\R^{n}\setminus\{0\}$ when $\Sigma=\S^{n-1}$. Our aim is to compute the best constant
$$
\delta_{n,\alpha}(\mathcal C_{\Sigma}):=\inf_{\scriptstyle u\in C^{2}_{c}(\overline{\mathcal C_{\Sigma}}\setminus\{0\})\atop\scriptstyle u\ne 0}\frac{\displaystyle\int_{\mathcal C_{\Sigma}}|x|^{\alpha}|\Delta u|^{2}dx}{\displaystyle\int_{\mathcal C_{\Sigma}}|x|^{\alpha-2}|\nabla u|^{2}dx}~\!.
$$
In fact this goal was already accomplished by Ghoussoub and Moradifam in \cite{GM} in the case of the whole space. 
However we provide alternative proofs which naturally adapt to handle with cone like domains. Even if this  generalization to cones seems to have a somehow artificial flavour,  in fact in our opinion it contains some deeper features. Firstly it allows us to consider the case of domains, even very regular, like the half-space, such that the singularity stays on the boundary. Moreover our results are stated in a fashion which makes clearer the expression of the best constant even in the case of the whole space. This fact is strongly related to the peculiar approach followed here. 

We also mention the papers \cite{AS} and \cite{Mus} dealing with a class of inequalities for radially symmetric functions on $\R^{n}$ in the non Hilbertian case, that is, involving the weighted $L^p$-norms of $-\Delta u$ and $\nabla u$, with $p>1$.

In order to state our main results we put
$$
\gamma_{n,\alpha}=\frac{(n-4+\alpha)(n-\alpha)}{4}~,\quad
  h_{n,\alpha}=\left(\frac{n-4+\alpha}{2}\right)^2.
$$
Given a domain $\Sigma$ in $\S^{n-1}$ with $\partial\Sigma\in C^{2}$, we denote by $\Lambda_{\Sigma}$ the spectrum of the Laplace-Beltrami operator on $\Sigma$ with null boundary conditions and by $\lambda_{\Sigma}$ the first eigenvalue. Notice that $\lambda_{\Sigma}>0$ apart from the case $\Sigma=\S^{n-1}$.

\begin{thm} 
\label{T:main}
Let $n\ge 2$ and let $\Sigma$ be a domain in $\S^{n-1}$ with $\partial\Sigma\in C^{2}$. Assume $\alpha\neq 4-n$. Then the following facts hold.
\begin{itemize}
\item[\textbf{(i)}] 
$\delta_{n,\alpha}(\mathcal C_{\Sigma})>0$ if and only if $-\gamma_{n,\alpha}\not\in\Lambda_{\Sigma}$. Moreover
\begin{equation}
\label{eq:above}
\delta_{n,\alpha}(\mathcal C_{\Sigma})\le M_{n,\alpha}(\Sigma):=
\min_{\lambda\in\Lambda_{\Sigma}}\frac{\left(\ga+\lambda\right)^2}{\ha+\lambda}~\!.
\end{equation}

\item[\textbf{(ii)}]
If $\gamma_{n,\alpha}-2h_{n,\alpha}\le\lambda_{\Sigma}$ then $\delta_{n,\alpha}(\mathcal C_{\Sigma})=M_{n,\alpha}(\Sigma)$.
\end{itemize}
\end{thm}
When $\Sigma=\S^{n-1}$ we can be more precise. First of all, as well as in \cite{GM}, we have the following sharp result for the best constant in the class of radially symmetric functions.

\begin{thm}
\label{T:radial}
Let $n\ge 2$ and $\alpha\in\R$. Then
$$
\delta_{n,\alpha}^{\rm rad}:=\inf_{\scriptstyle u\in C^{2}_{c}(\R^n\setminus\{0\})
\atop\scriptstyle u=u(|x|)~,~u\ne 0}\frac{\displaystyle\int_{\R^n}|x|^{\alpha}|\Delta u|^{2}dx}{\displaystyle\int_{\R^n}|x|^{\alpha-2}|\nabla u|^{2}dx}=\left(\frac{n-\alpha}{2}\right)^2~.
$$
\end{thm}
Theorem \ref{T:radial} will be proved in Section \ref{SS:radial}. 
When we allow $u$ to be any function in $C^{2}_{c}(\R^{n}\setminus\{0\})$ non necessarily radial we can estimate the best constant with the aid of Theorem \ref{T:main} and using the explicit knowledge of the spectrum of the Laplace-Beltrami operator on the sphere:
$$
\Lambda_{\S^{n-1}}=\{k(n-2+k)~|~k\in\mathbb{N}\cup\{0\}\}.
$$ 
In particular $\lambda_{\S^{n-1}}=0$. To simplify the notation, we write $\delta_{n,\alpha}$ instead of $\delta_{n,\alpha}(\R^{n}\setminus\{0\})$, and $M_{n,\alpha}$ instead of $M_{n,\alpha}(\S^{n-1})$. The results stated in the next theorems are already known (see \cite{GM}) but we prove them in a different way. 

\begin{thm} 
\label{T:mainRN}
Let $n\ge 2$ and assume $\alpha\neq 4-n$. 
\begin{itemize}
\item[\textbf{(i)}]
If $n=2$ then $\delta_{2,\alpha}=M_{2,\alpha}$ for any $\alpha\in\R$.
\item[\textbf{(ii)}]
If $n\ge 3$ then there exists $\alpha^*\in[4-n,2)$ such that $\delta_{n,\alpha}=M_{n,\alpha}$ for any $\alpha\notin[4-n,\alpha^*)$. 
\item[\textbf{(iii)}]
If $n\ge 3$ and $\alpha^*<\alpha<n$ then $\da=\darad$.
\end{itemize}
\end{thm}
In the ``\emph{critical case}'' $\alpha=4-n$ a very singular phenomenon can be observed.

\begin{thm} 
\label{T:4-n}
If $\alpha=4-n$ then
$$
\delta_{n,4-n}=\min\left\{(n-2)^2~,~n-1\right\}.
$$
In particular, 
\smallskip

\centerline{
$\delta_{n,4-n}>0$ for any $n\ge 3$ and
$\delta_{n,4-n}=n-1<\delta_{n,4-n}^{\rm rad}$ for any $n\ge 4$.}
\end{thm}
It should be emphasized the fact that the function $\alpha\mapsto\da$ is not continuous at $\alpha=4-n$, unless $n=2$.
Let us make some remarks about the above results in the meaningful case $\alpha=0$. First notice that in two dimensions $\delta_{2,0}=0<\delta_{2,0}^{\rm rad}=1$. In dimension $n=3$ the best constant $\delta_{3,0}$, already known according to the paper \cite{GM} can be computed by means of the formula for $M_{3,0}$ and yields:
$$
\int_{\R^3}|\Delta u|^2~dx\ge \frac{25}{36}\int_{\R^3}|x|^{-2}|\nabla u|^2~dx
\quad\textrm{for any $u\in C^2_c(\R^3\setminus\{0\})$.}
$$
Notice that $\delta_{3,0}^{\rm rad}=9/4$ is larger than the best constant on the whole space and breaking symmetry occurs. A similar phenomenon appears in the {\em critical dimension} $n=4$. Indeed
$\delta_{4,0}^{\rm rad}=4$, while from Theorem \ref{T:4-n} it follows that $3$ is the best constant in the inequality
$$
\int_{\R^4}|\Delta u|^2~dx\ge 3\int_{\R^4}|x|^{-2}|\nabla u|^2~dx
\quad\textit{for any $u\in C^2_c(\R^4\setminus\{0\})$.}
$$

To handle higher dimensions we estimate
\begin{equation}
\label{eq:alphastar}
\alpha^*<\frac{1}{3}\left(n+4-2\sqrt{n^2-n+1}\right).
\end{equation}
Notice that $\alpha^*<0$ if $n\ge 5$. A standard density result can
be used to infer the next corollary.

\begin{cor}
Assune $n\ge 5$. Then
$$
\irn|\Delta u|^2~dx\ge \frac{n^2}{4}\irn|x|^{-2}|\nabla u|^2~dx
\quad\textit{for any $u\in \mathcal D^{2,2}(\R^n)$,}
$$
and $n^2/4$ is the best constant.
\end{cor}

\section{Proofs}
The only tools we use are the Cauchy-Schwarz inequality, integration by parts, the variational characterization of the eigenvalues, and the Emden-Fowler transform $u\mapsto w=Tu$, that is defined via
$$
u(x)  =  |x|^{\frac{4-n-\alpha}{2}}~\!w\left(-\log|x|,\frac{x}{|x|}\right).
$$
Such a transform $T$ maps functions $u\colon{\R^n}\setminus\{0\}\to\R$ into functions $w=w(s,\sigma)$ on the cylinder $\R\times\S^{n-1}$. 
More generally, given a domain $\Sigma$ in $\S^{n-1}$, let us denote 
$$
\mathcal{Z}_{\Sigma}:=\R\times\Sigma
$$ 
the corresponding cylinder. We point out that $w\in C^{2}_{c}(\overline{\mathcal Z_{\Sigma}})$ as $u\in C^{2}_{c}(\overline{\mathcal C_{\Sigma}}\setminus\{0\})$. Moreover, by direct computation (see for instance \cite{CM1}), it can be proved that
$$
\da(\mathcal C_{\Sigma})=\inf_{\scriptstyle w\in C^{2}_{c}(\overline{\mathcal Z_{\Sigma}})\atop\scriptstyle w\ne 0}\frac
{\displaystyle \int_{\mathcal{Z}_{\Sigma}}\!\left
 |\Delta_\sigma w+w_{ss}+(\alpha-2)w_s-\ga w\right|^2 dsd\sigma}
{\displaystyle\int_{\mathcal Z_{\Sigma}}\!\left(|\nabla_\sigma w|^2+|w_s|^2\right) dsd\sigma+\ha \int_{\mathcal Z_{\Sigma}}|w|^2 dsd\sigma}~.
$$
Here and in the rest of the paper we denote by $-\Delta_\sigma$, $\nabla_{\!\sigma}$ the Laplace-Beltrami operator and the gradient on $\S^{n-1}$, respectively, while $w_s$ is the derivative of $w$ with respect to $s\in\R$.

\subsection{Some notation and technical lemmas}
For every eigenvalue $\lambda\in\Lambda_{\Sigma}$ let 
$$ 
Y_{\lambda}:=\{g\varphi~|~g\in C^{2}_{c}(\R),~\varphi~\textrm{eigenfunction corresponding to }\lambda\}.
$$
Notice that $Y_{\lambda}\subset C^{2}_{c}(\overline{\mathcal Z_{\Sigma}})$. Moreover set
$$
V_{\lambda}:=\{v\in C^{2}_{c}(\overline{\mathcal Z_{\Sigma}})~|~\int_{\mathcal Z_{\Sigma}}vw\!~dsd\sigma=0~\forall w\in Y_{\lambda'},\forall\lambda'\in\Lambda_{\Sigma},~\lambda'<\lambda\}.
$$
Hence $V_{\lambda}\supset Y_{\lambda}$ and $V_{\lambda}=C^{2}_{c}(\overline{\mathcal Z_{\Sigma}})$ when $\lambda=\lambda_{\Sigma}$. 
\begin{lem}\label{L:1}
For every $\lambda\in\Lambda_{\Sigma}$ one has that 
$$
\inf_{\scriptstyle w\in V_{\lambda}\atop\scriptstyle w\ne 0}\frac{\displaystyle\int_{\mathcal{Z}_{\Sigma}}\left(|\nabla_{\sigma}w|^{2}+|w_{s}|^{2}\right)~\!dsd\sigma}{\displaystyle\int_{\mathcal{Z}_{\Sigma}}|w|^{2}~\!dsd\sigma}\ge\lambda.
$$
\end{lem}

\proof 
Clearly 
$$
\inf_{\scriptstyle w\in V_{\lambda}\atop\scriptstyle w\ne 0}\frac{\displaystyle\int_{\mathcal{Z}_{\Sigma}}\left(|\nabla_{\sigma}w|^{2}+|w_{s}|^{2}\right)~\!dsd\sigma}{\displaystyle\int_{\mathcal{Z}_{\Sigma}}|w|^{2}~\!dsd\sigma}=
\inf_{\scriptstyle w\in V_{\lambda}\atop\scriptstyle w\ne 0}\frac{\displaystyle\int_{\mathcal{Z}_{\Sigma}}|\nabla_{\sigma}w|^{2}~\!dsd\sigma}{\displaystyle\int_{\mathcal{Z}_{\Sigma}}|w|^{2}~\!dsd\sigma}~\!.
$$
Then the conclusion follows from the fact that every mapping $w\in V_{\lambda}$ is orthogonal to $Y_{\lambda'}$ for any eigenvalue $\lambda'<\lambda$ and from the variational charachetrization of the eigenvalues.
\qed
\medskip

\noindent
For $A,B,C\in\R$ set
\begin{equation}\label{eq:ND}
\begin{array}{c}
\displaystyle N_{A,B}(w)=\int_{\mathcal{Z}_{\Sigma}}|\Delta_{\sigma}w+w_{ss}+Aw_{s}-Bw|^{2}~\!dsd\sigma\\
\displaystyle D_{C}(w)=\int_{\mathcal{Z}_{\Sigma}}\left(|\nabla_{\sigma}w|^{2}+|w_{s}|^{2}\right)~\!dsd\sigma+C\int_{\mathcal{Z}_{\Sigma}}|w|^{2}~\!dsd\sigma~\!.
\end{array}
\end{equation}
Moreover, for $\lambda\in\Lambda_{\Sigma}$, set
$$
M_{\lambda}(A,B,C)=\inf_{\scriptstyle w\in V_{\lambda}\atop\scriptstyle w\ne 0}\frac{N_{A,B}(w)}{D_{C}(w)}\quad\textrm{and}\quad\widetilde{M}_{\lambda}(A,B,C)=\inf_{\scriptstyle w\in Y_{\lambda}\atop\scriptstyle w\ne 0}\frac{N_{A,B}(w)}{D_{C}(w)}~\!.
$$

\begin{lem}\label{L:V}
For every $\lambda\in\Lambda_{\Sigma}$, if $0<B+\lambda\le 2(C+\lambda)$, then 
$$
M_{\lambda}(A,B,C)\ge\frac{(B+\lambda)^{2}}{C+\lambda}~\!.
$$
\end{lem}

\proof
For every $w\in V_{\lambda}$, integrating by parts and using Cauchy-Schwarz inequality, we obtain
\begin{equation*}
\begin{split}
D_{B}(w)&=-\int_{\mathcal{Z}_{\Sigma}}w\left(\Delta_{\sigma}w+w_{ss}+Aw_{s}-Bw\right)~\!dsd\sigma\\
&\le\left(\int_{\mathcal{Z}_{\Sigma}}|w|^{2}~\!dsd\sigma\right)^{\frac{1}{2}}\left(N_{A,B}(w)\right)^{\frac{1}{2}}.
\end{split}
\end{equation*}
Then for $w\in V_{\lambda}\setminus\{0\}$ we have
$$
\frac{N_{A,B}(w)}{D_{C}(w)}\ge\frac{\left(R(w)+B\right)}{R(w)+C}\quad\textrm{where}\quad R(w)=\frac{\displaystyle\int_{\mathcal{Z}_{\Sigma}}\left(|\nabla_{\sigma}w|^{2}+|w_{s}|^{2}\right)~\!dsd\sigma}{\displaystyle\int_{\mathcal{Z}_{\Sigma}}|w|^{2}~\!dsd\sigma}.
$$
Therefore, using Lemma \ref{L:1}, we infer that
$$
M_{\lambda}(A,B,C)\ge\inf_{r\ge\lambda}\frac{(B+r)^{2}}{C+r}=\frac{(B+\lambda)^{2}}{C+\lambda}
$$
where the last equality can be obtained by elementary calculus using the assumptions on $B$ and $C$. 
\qed

\begin{lem}\label{L:Y}
For every $\lambda\in\Lambda_{\Sigma}$, if $A^{2}+2(B+\lambda)>(B+\lambda)^{2}/(C+\lambda)$ and $C+\lambda>0$, then 
$$
\widetilde{M}_{\lambda}(A,B,C)\ge\frac{(B+\lambda)^{2}}{C+\lambda}~\!.
$$
\end{lem}

\proof
Using the definition of $Y_{\lambda}$ we obtain that
$$
\widetilde{M}_{\lambda}(A,B,C)=\inf_{\scriptstyle g\in C^{2}_{c}(\R)\atop\scriptstyle g\ne 0}\frac{\displaystyle\int_{-\infty}^{\infty}|g''+Ag'-(B+\lambda)g|^{2}~\!ds}{\displaystyle\int_{-\infty}^{\infty}\left(|g'|^{2}+(C+\lambda)|g|^{2}\right)~\!ds}~\!.
$$
To simplify notation, we can assume that
$$
\int_{-\infty}^{\infty}|g|^{2}~\!ds=1.
$$
Integration by parts yields
\begin{equation*}
\begin{split}
\int_{-\infty}^{\infty}|g''+Ag'-(B+\lambda)g|^{2}~\!&ds=
\int_{-\infty}^{\infty}|g''|^{2}~\!ds\\
&+(A^{2}+2(B+\lambda))\int_{-\infty}^{\infty}|g'|^{2}~\!ds+(B+\lambda)^{2}.
\end{split}
\end{equation*}
Moreover by Cauchy-Schwarz and Young inequality we estimate
$$
\int_{-\infty}^{\infty}|g'|^{2}~\!ds=-\int_{-\infty}^{\infty}g''g~\!ds\le\left(\int_{-\infty}^{\infty}|g''|^{2}~\!ds\right)^{\frac{1}{2}}\le\frac{\varepsilon}{2}+\frac
{1}{2\varepsilon}\int_{-\infty}^{\infty}|g''|^{2}~\!ds~\!.
$$
Then
\begin{equation*}
\begin{split}
&\frac{\displaystyle\int_{-\infty}^{\infty}|g''+Ag'-(B+\lambda)g|^{2}~\!ds}{\displaystyle\int_{-\infty}^{\infty}\left(|g'|^{2}+(C+\lambda)|g|^{2}\right)~\!ds}\\
&\qquad\qquad\qquad\ge
\frac{(A^{2}+2(B+\lambda+\varepsilon))\displaystyle\int_{-\infty}^{\infty}|g'|^{2}~\!ds+(B+\lambda)^{2}-\varepsilon^{2}}{\displaystyle\int_{-\infty}^{\infty}|g'|^{2}~\!ds+C+\lambda}
\end{split}
\end{equation*}
and consequently
$$
\widetilde{M}_{\lambda}(A,B,C)\ge(A^{2}+2(B+\lambda+\varepsilon))\inf_{t\ge 0}\frac{B_{\varepsilon}+t}{C+\lambda+t}
$$
where
$$
B_{\varepsilon}=\frac{(B+\lambda)^{2}-\varepsilon^{2}}{A^{2}+2(B+\lambda+\varepsilon)}.
$$
By the assumptions on $A$, $B$ and $C$ we have that $C+\lambda>B_{\varepsilon}>0$ for $\varepsilon>0$ small enough. Then, by elementary calculus, 
$$
\inf_{t\ge 0}\frac{B_{\varepsilon}+t}{C+\lambda+t}=\frac{B_{\varepsilon}}{C+\lambda}.
$$
Hence for $\varepsilon>0$ small enough
$$
\widetilde{M}_{\lambda}(A,B,C)\ge\frac{(B+\lambda)^{2}-\varepsilon^{2}}{C+\lambda}
$$
and letting $\varepsilon\to 0$ we get the conclusion.
\qed

\subsection{Proof of Theorem \ref{T:main}}

Fix $\alpha\in\R$ and $n\in\mathbb{N}$, $n\ge 2$. 
For every $w\in C^{2}_{c}(\overline{\mathcal Z_{\Sigma}})$ set
$$
\begin{array}{c}
\displaystyle N(w)=\int_{\mathcal{Z}_{\Sigma}}|\Delta_\sigma w+w_{ss}+(\alpha-2)w_s-\gamma_{n,\alpha}w|^{2}\!~dsd\sigma\\
\displaystyle D(w)=\int_{\mathcal Z_{\Sigma}}\left(|\nabla_\sigma w|^{2}+|w_s|^{2}\right)~\!dsd\sigma+h_{n,\alpha}\int_{\mathcal Z_{\Sigma}}|w|^{2}~\!dsd\sigma~\!.
\end{array}
$$
Notice that according to the notation (\ref{eq:ND}) we have that $N=N_{A,B}$ and $D=D_{C}$ with
\begin{equation}
\label{eq:ABC}
A=\alpha-2~\!,\quad B=\gamma_{n,\alpha}~\!,\quad C=h_{n,\alpha}~\!.
\end{equation}
\textbf{Proof of (i).} Since $\alpha\neq 4-n$, then $\ha>0$ and therefore the functional $D$ is the square of an equivalent Hilbertian norm on $H^1(\mathcal Z_{\Sigma})$. Assume that $-\ga\not\in\Lambda_{\Sigma}$. In this case, by the results in \cite{CM1}, the functional $N$ is the square of an equivalent Hilbertian norm on $H^2(\mathcal Z_{\Sigma})$. Therefore, since with the above notation 
$$
\da(\mathcal C_{\Sigma})=\inf\left\{\frac{N(w)}{D(w)}~\left|~w\in C^{2}_{c}(\overline{\mathcal Z_{\Sigma}}),~w\ne 0\right.\right\},
$$
as $H^2(\mathcal Z_{\Sigma})$ is continuously embedded into $H^1(\mathcal Z_{\Sigma})$, we obtain $\da(\mathcal C_{\Sigma})>0$. The fact that $\da(\mathcal C_{\Sigma})=0$ if $-\ga\in\Lambda_{\Sigma}$ is a consequence of (\ref{eq:above}). To check (\ref{eq:above}) we fix  $\lambda\in\Lambda_{\Sigma}$ and we estimate
\begin{equation}
\begin{split}
\nonumber
\da(\mathcal C_{\Sigma})&\le\inf_{\scriptstyle w\in Y_{\lambda}\atop\scriptstyle w\ne 0}\frac{N(w)}{D(w)}\\
\label{eq:above0}
&=\inf_{\scriptstyle g\in C^{2}_{c}(\R)\atop\scriptstyle g\ne 0}\frac{\displaystyle \int_{\R}\!\left |g''+(\alpha-2)g'-(\gamma_{n,\alpha}+\lambda)g\right|^2 ds}{\displaystyle\int_{\R}|g'|^2 ds+(\ha+\lambda)\int_{\R}|g|^2 ds}=\frac{(\gamma_{n,\alpha}+\lambda)^{2}}{h_{n,\alpha}+\lambda}.
\end{split}
\end{equation}
The last equality can be easily checked taking $g(s)=g_{0}(\varepsilon s)$ with $g_{0}\in C^{2}_{c}(\R)$ fixed, $g_{0}\ne 0$, and $\varepsilon>0$, and letting $\varepsilon\to 0$. Then (\ref{eq:above}) follows from the arbitrariness of $\lambda\in\Lambda_{\Sigma}$. 
\medskip

\noindent
\textbf{Proof of (ii).} It suffices to study the case $-\ga\not\in\Lambda_{\Sigma}$, since otherwise $\da(\mathcal C_{\Sigma})=\Ma(\Sigma)=0$. Let us distinguish the argument according that $-\gamma_{n,\alpha}$ stays below the spectrum or not.
\smallskip

\noindent
\textbf{Case $-\gamma_{n,\alpha}<\lambda_{\Sigma}$~.}
\smallskip

\noindent
Since $C^{2}_{c}(\overline{\mathcal Z_{\Sigma}})=V_{\lambda_{\Sigma}}$, we have that $\delta_{n,\alpha}(\mathcal C_{\Sigma})=M_{\lambda_{\Sigma}}(A,B,C)$ with $A$, $B$, and $C$ given as in (\ref{eq:ABC}). We apply Lemma \ref{L:V} with $\lambda=\lambda_{\Sigma}$. The condition $B+\lambda>0$ is fulfilled since we are dealing with the case $-\gamma_{n,\alpha}<\lambda_{\Sigma}$. The condition $B+\lambda\le 2(C+\lambda)$ is equivalent to say $\gamma_{n,\alpha}\le 2h_{n,\alpha}+\lambda_{\Sigma}$. Hence if $-\lambda_{\Sigma}<\gamma_{n,\alpha}\le 2h_{n,\alpha}+\lambda_{\Sigma}$ then
$$
\delta_{n,\alpha}(\mathcal C_{\Sigma})\ge\frac{(\gamma_{n,\alpha}+\lambda_{\Sigma})^{2}}{h_{n,\alpha}+\lambda_{\Sigma}}\ge M_{n,\alpha}(\Sigma).
$$
Hence, in this case, by (1), $\delta_{n,\alpha}(\mathcal C_{\Sigma})=M_{n,\alpha}(\Sigma)$. 
\smallskip

\noindent
\textbf{Case $-\gamma_{n,\alpha}>\lambda_{\Sigma}$~.}
\smallskip

\noindent
We can find two consecutive eigenvalues $\lambda_{k-1}$ and $\lambda_{k}$ such that
$$
\lambda_{k-1}<-\gamma_{n,\alpha}<\lambda_{k}~\!.
$$
Any $w\in C^{2}_{c}(\overline{\mathcal Z_{\Sigma}})$ can be written according to the following decomposition
$$
w=v_{1}+...+v_{k}
$$
with $v_{j}\in Y_{\lambda_{j}}$ for $j=1,...,k-1$, and $v_{k}\in V_{\lambda_{k}}$. One easily checks that
$$
\frac{N(w)}{D(w)}=\sum_{j=1}^{k}\theta_{j}\frac{N(v_{j})}{D(v_{j})}\quad\textrm{where }\theta_{j}=\frac{D(v_{j})}{D(w)}.
$$
Since $\theta_{j}\ge 0$ for all $j=1,...,k$ and $\theta_{1}+...+\theta_{k}=1$, we have that
\begin{equation}
\label{eq:estimate-w}
\frac{N(w)}{D(w)}\ge\min_{j=1,...,k}\frac{N(v_{j})}{D(v_{j})}.
\end{equation}
We estimate $N(v_{j})/D(v_{j})$ for $j=1,...,k-1$ by means of Lemma \ref{L:Y} with $\lambda=\lambda_{j}$ and $A$, $B$, and $C$ as in (\ref{eq:ABC}). The condition $C+\lambda>0$ is fulfilled as $\lambda_{j}\ge 0$ and $h_{n,\alpha}>0$ since, by hypothesis, $\alpha\ne 4-n$. The condition $A^{2}+2(B+\lambda)>(B+\lambda)^{2}/(C+\lambda)$ can be checked by considering the function
$$
\Phi(t)=\left(2t+\frac{(n-2)^{2}}{2}+\frac{(\alpha-2)^{2}}{2}\right)\left(t+h_{n,\alpha}\right)-\left(t+\gamma_{n,\alpha}\right)^{2}.
$$
One has that $\Phi(0)=h_{n,\alpha}^{2}>0$ and $\Phi'(0)=2h_{n,\alpha}+(\alpha-2)^{2}>0$. Then $\Phi(t)>0$ for all $t\ge 0$. In particular $\Phi(\lambda_{j})>0$ and 
$$
A^{2}+2(B+\lambda_{j})-\frac{(B+\lambda_{j})^{2}}{C+\lambda_{j}}=\frac{\Phi(\lambda_{j})}{h_{n,\alpha}+\lambda_{j}}>0.
$$
Hence Lemma \ref{L:Y} applies and yields
\begin{equation}
\label{eq:estimate-vj}
\frac{N(v_{j})}{D(v_{j})}\ge\frac{\left(\gamma_{n,\alpha}+\lambda_{j}\right)^{2}}{h_{n,\alpha}+\lambda_{j}}\ge M_{n,\alpha}(\Sigma)\quad\forall j=1,...,k-1.
\end{equation}
In order to estimate $N(v_{k})/D(v_{k})$ we apply Lemma \ref{L:V} with $\lambda=\lambda_{k}$  and $A$, $B$, and $C$ as in (\ref{eq:ABC}). The condition $B+\lambda>0$ is satisfied since $-\gamma_{N,\alpha}<\lambda_{k}$. The other condition $B+\lambda\le 2(C+\lambda)$ is also fulfilled since 
$$
2(C+\lambda)-(B+\lambda)=2h_{n,\alpha}-\gamma_{n,\alpha}+\lambda_{k}>2h_{n,\alpha}-\gamma_{n,\alpha}+\lambda_{\Sigma}>0
$$
by the assumption made in \textbf{(ii)}. Therefore Lemma \ref{L:V} applies and thus
\begin{equation}
\label{eq:estimate-vk}
\frac{N(v_{k})}{D(v_{k})}\ge\frac{\left(\gamma_{n,\alpha}+\lambda_{k}\right)^{2}}{h_{n,\alpha}+\lambda_{k}}\ge M_{n,\alpha}.
\end{equation}
In conclusion by (\ref{eq:estimate-w})--(\ref{eq:estimate-vk}) and by the arbitrariness of $w\in C^{2}_{c}(\overline{\mathcal Z_{\Sigma}})$ one concludes as in the first case.
\qed

\subsection{Proof of Theorem \ref{T:radial}}
\label{SS:radial}
For a fixed radial function $u\in C^2_c(\R^n\setminus\{0\})$ we introduce the radially symmetric function
$$
v(x)=|x|^{\frac{2-n-\alpha}{2}}u_r(x),
$$
where $u_r$ is the radial derivative of $u$. Then
\begin{eqnarray*}
\irn|x|^\alpha|\Delta u|^2~\!dx &=& 
\irn|x|^\alpha\left|\frac{n-\alpha}{2}|x|^{-1}u_r+|x|^{\frac{2-n-\alpha}{2}}v_r\right|~\!dx\\
&=&\left(\frac{n-\alpha}{2}\right)^2\irn|x|^{\alpha-2}|\nabla u|^2~\!dx+\irn|x|^{2-n}|\nabla v|^2~\!dx,
\end{eqnarray*}
since the double product vanishes:
$$
\irn|x|^{\frac{\alpha-n}{2}}v_r u_r~\!dx=\irn|x|^{1-n}vv_r~\!dx=c\int_0^\infty(v^2)_r~\!dr=0.
$$
The conclusion is immediate.
\qed

\subsection{Proof of Theorem \ref{T:mainRN}}
We apply Theorem \ref{T:main} considering that we deal with the case $\lambda_{\Sigma}=0$. 
\smallskip

\noindent
\textbf{Proof of (i).}
If $n=2$ and $\alpha\ne 2$ then $\gamma_{n,\alpha}=-(\alpha-2)^{2}/4<0$. The condition 
\begin{equation}
\label{eq:cond}
\gamma_{n,\alpha}-2h_{n,\alpha}\le\lambda_{\Sigma}
\end{equation} 
holds true for every $\alpha\ne 2$ and thus one can conclude. \smallskip

\noindent
\textbf{Proof of (ii).}
Consider now the case $n\ge 3$. Suppose $\gamma_{n,\alpha}>0$ i.e. $\alpha\in(4-n,n)$. In this case the condition (\ref{eq:cond}) holds true if and only if $\alpha\ge(n-8)/3$. When $\gamma_{n,\alpha}<0$, the condition (\ref{eq:cond}) always holds true. Hence \textbf{(ii)} is proved with $\alpha^{*}<(n-8)/3$.
\smallskip

\noindent
\textbf{Proof of (iii).}
If $n\ge 3$ and $\alpha\in(\alpha^{*},n)$ then $\gamma_{n,\alpha}>0$ and the mapping $t\mapsto(\gamma_{n,\alpha}+t)^{2}/(h_{n,\alpha}+t)$ is increasing in $[0,\infty)$. Hence $M_{n,\alpha}=\gamma_{n,\alpha}^{2}/h_{n,\alpha}=\delta_{n,\alpha}^{\textrm{rad}}$, by Theorem \ref{T:radial}.
\qed

\subsection{Proof of Theorem \ref{T:4-n}} 
\label{SS:1.2}
First notice that $\delta_{n,4-n}\le\delta_{n,4-n}^{\rm rad}=(n-2)^2$
by Theorem \ref{T:radial}. Now we prove that
$\delta_{n,4-n}\le n-1$. Notice that $\gamma_{n,4-n}=h_{n,4-n}=0$. We estimate $\delta_{n,4-n}$ with a family of mappings $w(s,\sigma)=g(\varepsilon s)\varphi(\sigma)$ where $g\in C^2_c(\R)$ is any nontrivial fixed function, $\varepsilon>0$ and $\varphi$ is an eigenfunction for $-\Delta_{\sigma}$ on $\S^{n-1}$ relative to the first positive eigenvalue $(n-1)$. In this way we obtain
$$
\delta_{n,4-n}\le\frac{\displaystyle \int_{\R}\!\left|\varepsilon^2 g''+(\alpha-2)\varepsilon g'-(n-1) g\right|^2 ds}{\displaystyle \varepsilon^2\int_{\R}|g'|^2~\!ds+(n-1)\int_\R|g|^2~\!ds}~\!.
$$
Then, passing to the limit as $\varepsilon\to 0$, we conclude that
$\delta_{n,4-n}\le n-1$. Thus $\delta_{n,4-n}\le \min\left\{(n-2)^2~,~n-1\right\}$. To prove the opposite inequality we argue by contradiction. We assume that there exists $w\in C^2_c(\R\times\S^{n-1})$, $w\neq 0$, such that
\begin{equation}
\label{eq:contradiction}
\min\left\{(n-2)^2~,~n-1\right\}>\frac{N(w)}{D(w)},
\end{equation}
where $N(w)$ and $D(w)$ are as in the proof of
Theorem \ref{T:main}.
We can write $w$ as
$w(s,\sigma)=g(s)+v(s,\sigma)$, where
$$
\int_{\S^{n-1}}v(s,\sigma)d\sigma=0
$$
for any $s\in\R$. Notice that $v\neq 0$, otherwise (\ref{eq:contradiction})
would contradict Theorem \ref{T:radial}.
Thus 
\begin{equation}
\label{eq:Poincare2}
\xi_v:=\frac{\displaystyle\int_{\mathcal Z}|\nabla_{\!\sigma} v|^2~\!dsd\sigma}
{\displaystyle\int_{\mathcal Z}|v|^2~\!dsd\sigma}\ge n-1.
\end{equation}
Clearly, $N(g)\ge(n-2)^2 D(g)$ by Theorem \ref{T:radial}. Arguing as in the proof of Theorem
\ref{T:main} and using (\ref{eq:Poincare2}) we can estimate
${N(v)}\ge (n-1){D(v)}$.
Therefore
\begin{eqnarray*}
\min\left\{(n-2)^2~,~n-1\right\}&>&\frac{N(w)}{D(w)}=
\frac{N(g)+N(v)}{D(g)+D(v)}\\
&\ge&
\frac{(n-2)^2 D(g)+(n-1)D(v)}{D(g)+D(v)},
\end{eqnarray*}
that readily leads to a contradiction. Thus equality holds and the theorem
is completely proved.
\qed

\subsection{Proof of (\ref{eq:alphastar})}
Let $\alpha\in(4-n,2)$ such that
$\darad>\da$. Thus there exist
$g\in C^2_c(\R)$ and $v\in C^2_c(\mathcal Z)$
such that $v(s,\cdot)$ has zero mean value on the sphere for any 
$s\in\R$, and such that 
$$
\darad>\frac{N(g+v)}{D(g+v)},
$$
where $N(\cdot)$ and $D(\cdot)$ are defined as in the proof of Theorem \ref{T:main}.
In addition, it holds that $v\neq 0$ and that $v$ satisfies (\ref{eq:Poincare2}).
Clearly, $N(g)\ge\darad D(g)$. Arguing as in the proof of Theorem
\ref{T:main} we can estimate
$$
\frac{N(v)}{D(v)}\ge\frac{(\xi_v+\ga)^2}{\xi_v+\ha},
$$
where $\xi_v$ is defined in (\ref{eq:Poincare2}). Therefore
$$
\darad>\frac{N(g)+N(v)}{D(g)+D(v)}\ge
\frac{\darad D(g)+\displaystyle\frac{(\xi_v+\ga)^2}{\xi+\ha}~\!D(v)}{D(g)+D(v)}.
$$
Noticing that $\ha\darad=\ga^2$, we infer that
$\darad-2\ga> \xi_v\ge n-1$. Hence $\alpha<2$ satisfies
$$
3\alpha^2-2(n+4)\alpha-n^2+4n+4> 0~,
$$
that is, 
$$\alpha< \frac{1}{3}~\left(n+4-2\sqrt{n^2-n-1}\right)~\!.$$
Conversely, if
$$
\frac{1}{3}~\left(n+4-2\sqrt{n^2-n-1}\right)\le\alpha<n,
$$
then it necessarily holds that $\da=\darad$.
\qed

\label{References}

\end{document}